\documentclass[11pt,a4paper,leqno]{amsart}
\input{avv.sty}
\title{Reflections on Ramanujan's\linebreak Mathematical Gems}
\author{G. D. Anderson and M. Vuorinen}
\subjclass[2010]{Primary 33-02, 33B15, 33C05, 33E05.\linebreak  Secondary 30C62.}

\begin{document}


\begin{abstract}  The authors provide a survey of certain aspects of their 
joint work with the late M. K. Vamanamurthy.
Most of the results are simple to state and deal with special 
functions, a topic of research where S. Ramanujan's contributions are 
well-known landmarks.  The comprehensive bibliography includes 
references to the latest contributions to this field.
\end{abstract}

\maketitle
\markboth{\textsc{G. D. ANDERSON AND
M. VUORINEN}}{\textsc{Reflections on Ramanujan's Mathematical Gems}}



\section{Introduction}

\subsection{Ramanujan's life}

Srinivasa Ramanujan (1887--1920), native to India, was an extraordinary mathematician. A child prodigy, he was largely self-educated. When he was 16, he found an 1856 book by G. S. Carr \cite{C}, that listed theorems and formulas and some short proofs.  Using this tutorial textbook, packed with facts from advanced calculus, geometry, and classical analysis, as a guide, Ramanujan taught himself mathematics, and by the age of 17 was engaged in deep mathematical research, studying Bernoulli numbers and divergent series and calculating the Euler-Mascheroni constant to 15 decimal places.

After several unsuccessful attempts to have his work appreciated by other mathematicians, he wrote to G. H. Hardy, who recognized his genius and invited him to study with him at Cambridge. Ramanujan, a devout Brahmin, at first refused to travel to a foreign country, but relented when the family goddess Namagiri appeared to his mother in a dream commanding her not to prevent his departure.

Ramanujan was tutored by, and collaborated with, Hardy for almost five years beginning in 1914. They published seven joint papers in varous journals \cite{Ram2}. Their brilliant guest made a deep impression on Hardy and Littlewood, who compared him to Jacobi and Euler. Hardy considered that his most important mathematical achievement was the discovery of Ramanujan.

Ramanujan was in poor health and was hospitalized for a long time because of a non-diagnosed illness. He returned in poor health to India in 1919, and died soon after at the age of 32. A modern analysis of his medical records has indicated that he may have been suffering from a form of hepatitis.

\subsection{Ramanujan's mathematical heritage}

It is beyond our competence to evaluate the significance of Ramanujan's mathematical genius.  In the literature, he is sometimes mentioned in the company of other great mathematicians such as Gauss, Jacobi, and Euler.  Ramanujan recorded most of his work in notebooks containing thousands of mathematical formulas and results.  In a series of books, published by Springer-Verlag from 1985 through 1998 \cite{Be1}--\cite{Be5}, B. Berndt carefully analyzed these notebooks, giving proofs for the results that Ramanujan had stated without proof.  Other work by Ramanujan is contained in the so-called Lost Notebook \cite{Ram1} and in some loose papers.

Berndt's analysis and book-writing project has required an enormous amount of effort and scientific detective work, for which he has received grateful acknowledgment from the mathematical community.  In particular, in 1996 he received the Steele Prize for Mathematical Exposition, with a citation \cite{S96} that reads, in part, ``In an impressive scholarly accomplishment spread out over 20 years, Berndt has provided a readable and complete account of the notebooks, making them accessible to other mathematicians.  Ramanujan's enigmatic, unproved formulas are now readily available, together with context and explication, often after the most intense and clever research efforts on Berndt's part.''

During the past ninety years that have passed since Ramanujan's death, his influence on several areas of mathematics such as number theory, combinatorics, and mathematical analysis has been significant and continues to be so.

\subsection{AVV meets Ramanujan}

For about twenty-five years the present authors had an active collaboration with the late M. K. Vamanamurthy, who died in 2009. We referred to our research group as AVV, after the initials of our last names. Our joint research dealt with geometric function theory, more precisely quasiconformal mapping theory. An important aspect of our work dealt with conformal invariants, usually expressed in terms of special functions such as the Euler gamma function, hypergeometric functions, complete elliptic integrals, and elliptic functions. 

By a lucky chance we discovered the survey of Askey \cite{As} and Berndt's series of books on Ramanujan's notebooks \cite{Be1}--\cite{Be5}, where we found valuable pieces of information. What interested us most was Ramanujan's work on the gamma and hypergeometric functions and modular equations. We also had access to a preprint version of \cite{BeBG}, which dealt with some of Ramanujan's theories.

These results of Ramanujan seemed to fit nicely into our AVV research program, started in about 1984--in particular, the part that dealt with special functions, \cite{AVV1}--\cite{AVV5} and \cite{ABRVV}, \cite{PV1}, \cite{AVV6}, \cite{AQVV}. 

\subsection{Old and new research}

   We have written several surveys on our AVV research. In \cite{ AVV1}
   we provided an overview of the known results, along with several
   new ones, and formulated a long list of open problems, including
   problems about the gamma function.  Next, in \cite{AVV5}, we outlined
   some of our earlier results and suggested that known inequalities
   and identities for the complete elliptic integral
   $$
\mathcal{K}(r) = \frac{\pi}{2}  {}_2F_1(\frac{1}{2}, \frac{1}{2};1;r^2)
$$
   might hold for  ${}_2F_1(a,b;a+b;r^2)$ with $(a,b)$ close to $(\frac{1}{2}, \frac{1}{2})\,,$ where ${}_2F_1(a,b;a+b;r^2)$ stands for the Gaussian
   hypergeometric function (see \cite{AS}, \cite{OLBC} and (\ref{eq:kolmas}) below).

It turned out that some of these
   ideas bore fruit, perhaps more than we had anticipated, in the following years, and the results are
   surveyed in \cite{AVV7} and \cite{AVV8}. Several of our research topics that were inspired by Ramanujan's work are the Euler-Mascheroni constant, the Euler gamma function, volumes of unit balls in euclidean n-space, approximation of the Gaussian hypergeometric function ${}_{2}F_{1}$, approximation of the perimeter of an ellipse, and the study of generalized modular equations.  Some of the most recent results, motivated by AVV work, include the papers authored by X.~Zhang, G. Wang, V. Heikkala, \'A. Baricz, E. A. Karatsuba, H. Alzer, and others.

The surveys \cite{AVV1} and \cite{AVV5} were written before our book \cite{AVV6}, which summarizes most of our work, whereas the surveys \cite{AVV7} and \cite{AVV8} were written after the publication of \cite{AVV6}.  The purpose of the present survey is to give a modified version of \cite{AVV7}, with an attempt to provide an overview of the most recent work on this subject matter.

While we where working on the book \cite{AVV6}, we became acquainted with the work of S.-L. Qiu, who subsequently visited each of the AVV team members at our respective home universities.  He helped us in checking the early versions of the book manuscript, and our collaboration with him led to many significant co-authored results in joint papers (\cite{AnQ}, \cite{AQVV}, \cite{QVa}, \cite{QVu1}, \cite{QVu2}, \cite{QVu3}).  Since then, his students and co-authors have energetically investigated the problems left open by our work, and also contributed in other ways to this area of research (\cite{QF}, \cite{QH}, \cite{QM}, \cite{QZ}, \cite{WZC}, \cite{WZQC}, \cite{ZWC1,ZWC2}).

\section{Gamma function and Euler-Mascheroni constant}

Throughout this paper $\Gamma$ will denote Euler's gamma function, defined by
\begin{equation*}
\Gamma (z) = \int^{\infty}_0 e^{-t} t^{z-1}\,dt, \quad  \text{Re}\,z > 0,
\end{equation*}
and then continued analytically to the finite complex plane minus the set of
nonpositive integers. The recurrence formula $\Gamma(z+1)=z\,\Gamma(z)$
yields $\Gamma(n+1)=n!$ for any positive integer $n$.
We also use the fact that $\Gamma(1/2)=\sqrt{\pi}.$
The beta function is related to the gamma function by
$B(a,b)= \Gamma(a)\Gamma(b)/\Gamma(a+b)$.
The logarithmic derivative of the gamma function will be denoted,
as usual, by
\begin{equation*}
\Psi(z)\equiv\frac{d}{dz} \log \Gamma(z) = \frac{\Gamma'(z)}{\Gamma(z)}.
\end{equation*}

The Euler-Mascheroni constant $\gamma$ is defined as (see \cite{A1},
\cite{TY}, \cite{Yo})
\begin{equation*}
\gamma \equiv \lim_{n\to\infty} D_n =0.57721 56649\dots ; \
D_n \equiv \sum_{k=1}^n\frac{1}{k}-\log n.
\end{equation*}

\noindent The convergence of the sequence
$D_n$ to $\gamma$ is very slow (the speed of convergence is studied
by Alzer \cite{A1}).  D. W. DeTemple \cite{De} studied
a modified sequence which
converges faster and proved
\begin{equation*}
\frac{1}{24(n+1)^2} < R_n-\gamma < \frac{1}{24n^2},
\quad \text{where} \quad
R_n\equiv \sum_{k=1}^n\frac{1}{k} - \log\left(n+\frac{1}{2}\right).
\end{equation*}

Now let
$$
h(n)=R_n-\gamma,\ \ H(n)=n^2h(n),\ n\geqslant 1.
$$
Since $\Psi(n) =-\gamma-1/n+ \sum_{k=1}^n{1/k},$ we see that
$$
H(n)=(R_n-\gamma)n^2 = \left(\Psi(n)+\frac{1}{n} -\log\left(n+ \frac{1}{2}\right)\right)n^2.
$$
Some computer experiments led M. Vuorinen to conjecture that
$H(n)$ increases on the interval $[1,\infty)$
from $H(1) = -\gamma + 1 - \log (3/2) = 0.0173\ldots$ to
$1/24 = 0.0416\ldots$.  E. A. Karatsuba proved in \cite{K1}
that for all integers $n \geqslant 1, H(n) < H(n+1),$
by clever use of Stirling's formula
and Fourier series.  Moreover, using the relation
$\gamma = 1-\Gamma'(2)$ she obtained, for $ k \geqslant 1$,
$$
-c_k \leqslant
\gamma -1 + (\log k)
 \sum_{r=1}^{12k+1} d(k,r)
- \sum_{r=1}^{12k+1}\frac{d(k,r)}{r+1}
\leqslant c_k ,\\
$$
where
$$
c_k=\frac{2}{(12k)!} +2 k^2 e^{-k},\quad
d(k,r)= (-1)^{r-1} \frac{k^{r+1}}{(r-1)! (r+1)},
$$
giving exponential convergence.
Some computer experiments also seemed to indicate that
$((n+1)/n)^2 H(n)$ is a decreasing convex function.

\begin{Rama}
In ``The Lost Notebook and Other Unpublished Papers'' of Ramanujan \cite{Ram1},
appears the following record:
\begin{equation*}
``\Gamma (1+x) =\sqrt{\pi} \Big(\frac{x}{e}\Big)^x
\Big\{8x^3+4x^2+x+\frac{\theta_x}{30}\Big\}^{1/6},
\end{equation*}
where $\theta_x$ is a positive proper fraction
\begin{align*}
&\theta_0 =\frac{30}{\pi^3}=.9675\\
&\theta_{1/12} =.8071\quad \theta_{7/12}=.3058\\
&\theta_{2/12} =.6160\quad \theta_{8/12}=.3014\\
&\theta_{3/12} =.4867\quad \theta_{9/12}=.3041\\
&\theta_{4/12} =.4029\quad \theta_{10/12}=.3118\\
&\theta_{5/12} =.3509\quad \theta_{11/12}=.3227\\
&\theta_{6/12} =.3207\quad \theta_1=.3359\\
&\theta_{\infty }=1.\text{''}
\end{align*}
\end{Rama}

Of course, the values in the above table, except $\theta_{\infty}$,
are irrational and hence the decimals should be nonterminating as
well as nonrecurring.  The record stated above has been
the subject of intense investigations and is reviewed
in \cite[page 48, Question 754]{BCK}.  This note of Ramanujan
led the authors of \cite{AVV6} to make the following conjecture.

\begin{conj}\label{conj:eka}
Let
\begin{equation*}
G(x) = (e/x)^x\Gamma (1+x)/\sqrt{\pi}
\end{equation*}
and
\begin{equation*}
H(x)=G(x)^6-8x^3-4x^2-x = \frac{\theta_x}{30}.
\end{equation*}
Then $H$ is increasing from $(1,\infty)$ into $(1/100,1/30)$
\cite[p.\ 476]{AVV6}.
\end{conj}
\begin{Kara}
In a nice piece of work, E. A. Karatsuba \cite{K2} proved 
Conjecture \ref{conj:eka}. She did so by representing the function $H(x)$
as an integral for which she was able to  find an asymptotic
development. Her work also led to an interesting asymptotic formula
for the gamma function:
\begin{equation}\label{Kara}
\begin{split}
&\Gamma (x+1)=\sqrt{\pi }\Big(\frac{x}{e}\Big)^x \Big(8x^3+4x^2+x+
\frac{1}{30}-\frac{11}{240x} + \frac{79}{3360x^2}  +
\frac{3539}{201600x^3}\\
&-\frac{9511}{403200x^4}-\frac{10051}{716800x^5}+
\frac{47474887}{1277337600x^6}
+\frac{a_7}{x^7}+\cdots+\frac{a_n}{x^n} +\Delta_{n+1}(x)\Big)^{1/6},
\end{split}
\end{equation}
\medskip
where $\Delta_{n+1}(x)=O(\frac{1}{x^{n+1}})$, as $x \to \infty$,
and where each $a_k$ is given explicitly in terms of the Bernoulli numbers.
\end{Kara}

 G. Nemes has studied the Ramanujan-Karatsuba formula in (\ref{Kara}) and shown that it is better than some other well-known approximations for the gamma function \cite{N}.

The {Monotone l'H\^opital's Rule},
stated in the next paragraph, played an important
role in our work \cite{AVV4}--\cite{AVV5}.
The authors discovered this result in \cite{AVV4}, unaware
that it had been used earlier (without the name) as a technical tool
in differential geometry.  See \cite[p.\ 124, Lemma 3.1]{Cha}
or \cite[p.\ 14]{AQVV} for relevant remarks.

\begin{lem}\label{lem:eka}
For $-\infty < a < b < \infty$,
let  $g$  and  $h$  be real-valued functions that are
continuous on $[a,b]$ and differentiable on $(a,b)$, with $h' \ne 0$ on
$(a,b)$.  If  $g'/h'$ is strictly increasing
{\rm (}resp. decreasing{\rm )} on $(a,b)$, then the functions
\begin{equation*}
\frac{g(x)-g(a)}{h(x)-h(a)}\quad \text{and}\quad \frac{g(x)-g(b)}{h(x)-h(b)}
\end{equation*}
are also strictly increasing {\rm (}resp. decreasing{\rm )} on $(a,b)$.
\end{lem}

\begin{MonProp}
In \cite{AnQ} it is proved that the function
\begin{equation}\label{eq:eka}
f(x)\equiv \frac{\log \Gamma (x+1)}{ x\, \log x}
\end{equation}
is strictly
increasing from $(1,\infty )$ onto $(1-\gamma,1)$.
In particular, for $x\in (1,\infty )$,
\begin{equation}\label{eq:toka}
x^{(1-\gamma )x-1} < \Gamma (x) < x^{x-1}.
\end{equation}
The proof required the following two technical
lemmas, among others:
\end{MonProp}

\begin{lem}\label{lem:toka}
The function
\begin{equation*}
g(x)\equiv \sum_{n=1}^{\infty }
\frac{n-x}{(n+x)^3}
\end{equation*}
is positive for $x\in [1,4)$.
\end{lem}

\begin{lem}\label{lem:kolmas}
The function
\begin{equation}\label{eq:hx}
h(x)\equiv x^2\,\Psi'(1+x) -x\,\Psi (1+x) + \log \Gamma (1+x)
\end{equation}
is positive for all
$x\in [1,\infty )$.
\end{lem}

It was conjectured in \cite{AnQ} that the function $f$ in
\eqref{eq:eka} is concave on $(1,\infty )$.

\begin{HA}
Horst Alzer \cite{A1} has given an elegant proof of the
monotonicity of the function $f$ in \eqref{eq:eka}
by using the Monotone l'H\^opital's Rule \ref{lem:eka} and the convolution theorem for
Laplace transforms. In a later paper \cite{A3} he has improved
the estimates in \eqref{eq:toka} to
\begin{equation}\label{eq:2-10}
x^{\alpha(x-1) - \gamma} < \Gamma(x) < x^{\beta(x-1) - \gamma},
\quad x \in (0,1),
\end{equation}
where $\alpha \equiv 1 - \gamma = 0.42278\dots$,
$\beta \equiv \frac{1}{2}\left(\pi^2/6 - \gamma\right) =
0.53385\dots$ are best possible.
If $x \in (1, \infty)$, he also showed that \eqref{eq:2-10} holds
with best constants
$\alpha \equiv \frac{1}{2} \left(\pi^2/6 - \gamma\right)=0.53385\dots$,
$\beta \equiv 1$.
\end{HA}

\begin{EL}
Elbert and Laforgia  \cite{EL} have shown that the function
$g$ in Lemma \ref{lem:toka}
is positive for all $x>-1$.  They used this result to prove that
the function $h$ in Lemma \ref{lem:kolmas} is strictly decreasing
from $(-1,0]$ onto $[0,\infty)$ and strictly increasing from
$[0,\infty)$  onto $[0,\infty)$. They also showed that
$f'' < 0$ for $x>1$, thus proving the Anderson-Qiu
conjecture \cite{AnQ}, where $f$ is as in \eqref{eq:eka}.
\end{EL}

\begin{BP}
Berg and Pedersen \cite{BP} have shown that the function $f$ in
\eqref{eq:eka} is not only strictly increasing from $(0,\infty)$
onto $(0,1)$, but is even a (nonconstant) so-called \emph{Bernstein function}.
That is, $f > 0$ and $f'$ is completely monotonic, i.e.,
$f' > 0$, $f'' < 0$, $f'''> 0$, \dots.
In particular, the function $f$ is strictly
increasing and strictly concave on $(0,\infty )$.

In fact, they have proved the stronger result that $1/f$ is a
Stieltjes transform, that is, can be written in the form
$$ \frac{1}{f(x)} = c + \int^{\infty}_0 \frac{d \sigma(t)}{x+t}, \quad x > 0,$$
where the constant $c$ is non-negative and $\sigma$ is a non-negative
measure on $[0,\infty)$ satisfying
$$ \int^{\infty}_0 \frac{d \sigma(t)}{1 + t} < \infty.$$
In particular, for $1/f$ they have shown by using Stirling's formula
that $c = 1$.  Also they have obtained $d \sigma(t) = H(t) dt$,
where $H$ is the continuous density
$$ H(t) = \left\{ \begin{array}{ll}
t \dfrac{\log | \Gamma(1-t) | + (k-1) \log t}{(\log |\Gamma (1-t)|)^2 +
(k-1)^2 \pi^2}, &  t \in (k-1,k), k=1,2,\dots, \\[.4cm]
0 \hspace{1.9in}, & t = 1,2,\dots .  \end{array}\right. $$
Here log denotes the usual natural logarithm.  The density $H(t)$ tends
to $1/\gamma$ as $t$ tends to $0$, and $\sigma$ has no mass at $0$.
\end{BP}
\begin{rem}
In a series of papers I. Pinelis (see, e.g. \cite{P}) has advocated the
  use of the  Monotone l'H\^opital's Rule, Lemma \ref{lem:eka}. Probably partly because of
  his work,  during the past few years this result has found numerous applications to
  the study of special functions.
  In a forthcoming paper \cite{KVV}, a long list of papers is provided in which Jordan's inequality is refined.  Most of these refinements make use of Lemma \ref{lem:eka}.
\end{rem}

\medskip

\section{Volumes of balls}

Formulas for geometric objects, such as volumes of solids and arc lengths of
curves, often involve special functions.  For example, if $\Omega_n$
denotes the volume of the
unit ball $B^n=\{x:|x|< 1\}$ in $\R^n$, and if $\omega_{n-1}$ denotes
the $(n-1)$-dimensional surface area of the unit sphere $
S^{n-1} = \{x:|x|=1\}$, $n\geqslant 2$, then
\begin{equation*}
\Omega_n = \frac{\pi^{n/2}}{\Gamma ((n/2)+1)} \, ; \
{\omega_{n-1}} = n \Omega_n .
\end{equation*}

\noindent It is well known that for $n \geqslant 7$ both $\Omega_n$ and $\omega_n$
decrease to $0$ (cf. \cite[2.28]{AVV6}).
However, neither $\Omega_n$ nor $\omega_n$ is monotone for
$n$ on $[2,\infty)$. On the other hand,
$\Omega_n^{1/(n \log n)}$ decreases to $e^{-1/2}$ as $n \to \infty$
\cite[Lemma 2.40(2)]{AVV1}.

In 2000 H. Alzer \cite{A2} obtained the best possible constants
$a$, $b$, $A$, $B$, $\alpha$, $\beta$ such that
\begin{align*}
a\,\Omega^{\frac{n}{n+1}}_{n+1} &\leqslant \Omega_n \leqslant
b\,\Omega^{\frac{n}{n+1}}_{n+1},\\
\sqrt{\frac{n+A}{2 \pi}} &\leqslant \frac{\Omega_{n-1}}{\Omega_{n}}
\leqslant \sqrt{\frac{n+B}{2 \pi}},\\
\left(1 + \frac{1}{n} \right)^{\alpha} &\leqslant
\frac{\Omega^2_n}{\Omega_{n-1} \Omega_{n+1}} \leqslant
\left( 1 + \frac{1}{n} \right)^{\beta}
\end{align*}
for all integers $n \geqslant 1$. He showed that
$a = 2/\sqrt{\pi} = 1.12837\ldots$, $b = \sqrt{e} = 1.64872\dots$,
$A = 1/2$, $B= \pi/2 - 1 = 0.57079\dots$,
$\alpha = 2 - (\log \pi)/\log 2 = 0.34850\dots$, $\beta = 1/2$.
For some related results, see \cite{KlR}.

Recently H. Alzer \cite{A3} has obtained several sharp inequalities for $\Omega_{n}$. In particular, he showed that

$$
\frac{A}{\sqrt{n}}\leqslant(n+1)\frac{\Omega_{n+1}}{\Omega_{n}}-n\frac{\Omega_{n}}{\Omega_{n-1}}<\frac{B}{\sqrt{n}},\ {\rm for}\, \, n\geqslant 2,
$$
with the best possible constants $A=(4-\pi)\sqrt{2}=1.2139\ldots$ and $B=\frac{1}{2}\sqrt{2\pi}=1.2533\ldots$, refining and complementing work in \cite{KlR}.

The most recent studies dealing with the monotonicity properties of $\Omega_{n}$
 include the following papers: \cite{A4}, \cite{BP2}, \cite{M1}, \cite{QG}.

\section{Hypergeometric functions}\label{sect:hypfun}

Given complex numbers
$a$, $b$, and $c$ with $c\neq 0,-1,-2, \dots $,
the \emph{Gaussian hypergeometric function} is the analytic
continuation to the slit plane $\C\setminus [1,\infty)$ of
\begin{equation}\label{eq:kolmas}
F(a,b;c;z)\!= \!{}_2 F_1(a,b;c;z)\! \equiv\!
\sum_{n=0}^{\infty} \frac{(a,n)(b,n)}{(c,n)} \frac{z^n}{n!},
\quad |z|<1.
\end{equation}
Here $(a,0)=1$ for $a\neq 0$, and $(a,n)$
is the \emph{shifted factorial function}
\begin{equation*}
(a,n)\equiv a(a+1)(a+2) \cdots (a+n-1)
\end{equation*}
for $n=1,2,3,\ldots $.

The hypergeometric function $w = F(a,b;c;z)$ in \eqref{eq:kolmas}
has the simple differentiation formula
\begin{equation}\label{eq:neljas}
\frac{d}{dz}\,F(a,b;c;z)=\frac{ab}{c}\,F(a+1,b+1;c+1;z).
\end{equation}

The behavior of the hypergeometric function near $z=1$ in the
three cases $ a+b < c$,
$a+b = c$, and $ a+b > c, ~ a,b,c > 0$, is given by
\begin{equation}\label{eq:viides}
\begin{cases}
F(a,b;c;1) = \frac{\Gamma(c) \Gamma(c-a-b)}{\Gamma(c-a) \Gamma(c-b)},
\ a+b < c,\\
B(a,b)F(a,b;a+b;z)+\log(1-z) \\
\qquad \qquad \quad = R(a,b)+ O((1-z)\log(1-z)),\\
F(a,b;c;z) = (1-z)^{c -a -b} F(c-a,c-b;c;z),\ c < a+b,
\end{cases}
\end{equation}
where 
\begin{equation}\label{eq:RamConst}
R(a,b) = -2 \gamma - \Psi(a) - \Psi(b),\ \ R(a) \equiv R(a,1-a),\ \ R(\frac{1}{2}) = \log 16,
\end{equation}
and where log
denotes the principal branch of the complex logarithm.
The above asymptotic formula for the \emph{zero-balanced} case $a+b=c$
is due to Ramanujan (see \cite{As}, \cite{Be1}). This formula is
implied by \cite[15.3.10]{AS}.

The asymptotic formula \eqref{eq:viides} gives a precise description
of the behavior of
the function $F(a,b;a+b;z)$ near the logarithmic singularity $z=1$.
This singularity can be removed by an exponential change of variables,
and the transformed function will be nearly linear.
 
In  \cite{QF} it is shown that  $f(x)\equiv R(x)\sin(\pi x)$ is decreasing from $(0,1/2]$ onto $(\pi, \log 16]$, where  the  \emph{Ramanujan constant}  $R(x)$ is as defined in (\ref{eq:RamConst}).

\begin{thm}\label{thm:eka}
\cite{AQVV}
For  $a,b > 0$, let $k(x) = F(a,b;a+b;1-e^{-x})$, $x > 0$.
Then $k$ is an increasing and convex function with $k'((0, \infty)) =
(ab/(a+b)$, $\Gamma(a+b)/(\Gamma(a) \Gamma(b)))$.
\end{thm}

\begin{thm}\label{thm:toka}
\cite{AQVV}
Given $a,b > 0$, and  $a+b > c$, $d \equiv a+b-c$, the function
$\ell (x) = F(a,b;c;1-(1+x)^{-1/d})$, $ x > 0$, is increasing and
convex, with $\ell '((0,\infty)) = (ab/(cd)$,
$\Gamma(c)\Gamma(d)/(\Gamma(a)\Gamma(b)))$.
\end{thm}

\begin{GaussCont}
The six functions $F(a\pm 1, b;c;z)$, $F(a, b\pm 1; c;z)$,
$F(a, b;c \pm 1;z)$ are
said to be \emph{contiguous} to $F(a,b;c;z)$. Gauss discovered 15 relations
between $F(a,b;c;z)$ and pairs of its contiguous functions
\cite[15.2.10--15.2.27]{AS}, \cite[Section 33]{Rai2}.
If we apply these relations to the differentiation formula \eqref{eq:neljas}, we obtain
the following useful formulas.
\end{GaussCont}

\begin{thm}\label{thm:kolmas}
For $a,b, c > 0$, $z \in (0,1)$,  let $u = u(z) =
F(a-1,b;c;z)$, $v = v(z) = F(a,b;c;z)$,
$u_1 = u(1-z)$, $v_1 = v(1-z)$.  Then
\begin{align}
z\frac{du}{dz} &= (a -1) (v - u),\label{eq:kuudes}\\
z(1-z)\frac{dv}{dz} &= (c -a) u + (a-c+bz)  v,\label{eq:seits}
\end{align}
and
\begin{equation}\label{eq:kahd}
\frac{ab}{c} z(1-z) F(a+1,b+1;c+1;z) =(c-a) u
+(a-c+bz)v.
\end{equation}
Furthermore,
\begin{equation}\label{eq:kahd2}
 z(1\!-\!z) \dfrac{d}{dz}\bigl(uv_1 \!+\! u_1 v \!-\! vv_1 \bigr)\! =\!  (1\! -\! a\! -\! b) \bigl[(1\!-\!z)u v_1 \!-\! z u_1 v\! -\! (1\!-\!2z)v v_1) \bigr].
\end{equation}
\end{thm}

Formulas \eqref{eq:kuudes}-\eqref{eq:kahd} in Theorem \ref{thm:kolmas} are
well known. See, for example, \cite[2.5.8]{AAR}.  On the other hand, formula
\eqref{eq:kahd2}, which follows from \eqref{eq:kuudes}-\eqref{eq:seits}
is first proved in  \cite[3.13 (4)]{AQVV}.

Note that the formula
\begin{equation}\label{eq:yhd}
z(1-z) \frac{dF}{dz} = (c-b)F(a,b-1;c;z)+(b-c+az)F(a,b;c;z)
\end{equation}
follows from  \eqref{eq:seits} if we use
the symmetry property $F(a,b;c;z) = F(b,a;c;z)$.

\bigskip

\begin{cor}\label{cor:3.13}
With the notation of Theorem {\rm \ref{thm:kolmas}},
if $a \in (0,1),~ b= 1-a < c,$ then
$$ uv_1 + u_1 v - vv_1 = u(1) =
\frac{(\Gamma (c))^2} {\Gamma (c+ a-1) \Gamma(c-a+1)}.$$
\end{cor}

\medskip

\medskip

\section{Hypergeometric differential equation
}\label{sect:HDE}
\medskip

The function $F(a,b;c;z)$ satisfies
the hypergeometric differential equation
\begin{equation}\label{eq:kymm}
z(1-z)w''+[c-(a+b+1)z]w'-abw=0.
\end{equation}
Kummer discovered solutions of \eqref{eq:kymm} in various domains,
obtaining 24 in all; for a complete list of his solutions
see \cite[pp.\ 174, 175]{Rai2}.

\begin{lem}\label{lem:hde1}
{\rm (1)} If $2c=a+b+1$ then both  $F(a,b;c;z)$ and  $F(a,b;c;1-z)$
satisfy {\rm (\ref{eq:kymm})} in the lens-shaped region
$\{z:0<|z|<1,\ 0<|1-z|<1\}.$\linebreak
{\rm (2)} If $2c=a+b+1$ then both $F(a,b;c;z^2)$ and $F(a,b;c;1-z^2)$
satisfy the differential equation
\begin{equation}\label{eq:kymm2}
z(1-z^2)w''+[2c-1-(2a+2b+1)z^2]w'-4abz w=0
\end{equation}
in the common part of the disk $\{z:|z|<1\}$ and the lemniscate
$\{z:|1-z^2|<1\}$.
\end{lem}

       {\bf Proof.}
By Kummer (cf. \cite[pp. 174-177]{Rai2}), the functions
$F(a,b;c;z)$ and $F(a,b;a+b+1-c;1-z)$ are solutions of (\ref{eq:kymm})
in $\{z:0<|z|<1\}$ and $\{z:0<|1-z|<1\}$, respectively.
But $a+b+1-c = c$ under the stated hypotheses.  The result (2)
follows from result (1) by the chain rule.\ \ $\square$

\medskip

\begin{lem}\label{lem:hde2}
The function $F(a,b;c;\sqrt{1-z^2})$ satisfies the differential equation
$$
Z^3(1-Z)zw''-\{Z(1-Z)+[c-(a+b+1)Z]Zz^2\}w'-abz^3w=0,
$$
in the subregion of the right half-plane bounded by the lemniscate  $r^2 = 2\cos (2\vartheta )$, $-\pi /4 \leqslant \vartheta \leqslant \pi /4$, $z=re^{i\vartheta}$.  Here $Z=\sqrt{1-z^2}$, where the square root indicates the principal branch.
\end{lem}

{\bf Proof.}
From (4.1), the differential equation for $w=F(a,b;c;t)$ is given by
$$t(1-t)\frac{d^2w}{dt^2}+[c-(a+b+1)t]\frac{dw}{dt}-abw=0.$$
Now put $t=\sqrt{1-z^2}.$  Then
$$
\frac{dz}{dt}=-\frac{t}{z},\ \frac{dt}{dz}=-\frac{z}{t},\ \frac{d^2t}{dz^2}=-\frac{1}{t^3}
$$
and
$$
\frac{dw}{dt}=-\frac{t}{z}\frac{dw}{dz},\ \frac{d^2w}{dt^2}=\frac{t^2}{z^2}\frac{d^2w}{dz^2}-\frac{1}{z^3}\frac{dw}{dz}.
$$
So
$$
t(1-t)\Big[\frac{t^2}{z^2}w''-\frac{1}{z^3}w'\Big]+\Big[c-(a+b+1)t\Big]
\Big(-\frac{t}{z}\Big)w'-abw=0.
$$
Multiplying through by $z^3$ and replacing $t$ by $Z\equiv \sqrt{1-z^2}$ gives the result.$ \ \square$

\bigskip

If $w_1$ and $w_2$ are two solutions of a second order differential
equation, then their \emph{Wronskian} is defined to be
$W(w_1,w_2)\equiv w_1w_2'-w_2w_1'$.

\begin{lem}\label{lem:linind}
\cite[Lemma 3.2.6]{AAR}
If $w_1$ and $w_2$ are two linearly independent solutions of
\eqref{eq:kymm}, then
\begin{equation*}
W(z)=W(w_1,w_2)(z)=\frac{A}{z^c(1-z)^{a+b-c+1}},
\end{equation*}
where $A$ is a constant.
\end{lem}

(Note the misprint in \cite[(3.10)]{AAR}, where the coefficient
$x(1-x)$ is missing from the first term.)

\begin{lem}\label{lem:vv1}
If $2c=a+b+1$ then, in the notation of Theorem {\rm \ref{thm:kolmas}},
\begin{equation}\label{eq:duren}
 (c-a)(uv_1 + u_1v) + (a-1)vv_1 = A \cdot z^{1-c}(1-z)^{1-c}.
\end{equation}
\end{lem}

\medskip

For a proof see \cite{AVV7}.  Note that in the particular case $c=1, a=b= \frac{1}{2}$
the right side of (\ref{eq:duren}) is constant and the result is similar to
Corollary \ref{cor:3.13}. This particular case is
Legendre's Relation (\ref{eq:legendre}; an elegant proof of it was given by
Duren \cite{Du}.

\medskip

\begin{lem}\label{lem:vv11}
 If $a, b > 0, c \geqslant 1,$ and $2c = a+b+1,$
then the constant $A$ in Lemma {\rm \ref{lem:vv1}} is given by
$A = (\Gamma(c))^2 / (\Gamma(a) \Gamma(b)).$
In particular, if $c = 1$ then Lemma {\rm \ref{lem:vv1}} reduces
to Legendre's Relation {\rm (\ref{eq:legkaea})} for generalized elliptic integrals.
\end{lem}

For a detailed proof of this lemma we refer the reader to \cite{AVV7}.
\bigskip

For rational triples $(a,b,c)$ there are numerous cases
where the hypergeometric function $F(a,b;c;z)$ reduces to a simpler
function (see \cite{PBM}). Other important particular cases are
\emph{generalized elliptic integrals}, which we will now discuss.
For $a,r\in(0,1)$, the \emph{generalized elliptic integral
of the first kind} is given by
\begin{align*}
\mathcal{K}_a &= \mathcal{K}_a(r) = \frac{\pi}{2}
F(a,1-a;1;r^2)\\
&= (\sin\pi a)\int_0^{\pi /2}(\tan t)^{1-2a}(1-r^2\sin^2t)^{-a}\,
dt,\\
\mathcal{K}_a' &=\mathcal{K}_a'(r)=\mathcal{K}_a(r').
\end{align*}
We also define
\begin{equation*}
\mu_a(r)=\frac{\pi}{2\sin (\pi a)}
\frac{\mathcal{K}_a'(r)}{\mathcal{K}_a(r)},\quad r'=\sqrt{1-r^2}.
\end{equation*}
For $\mu_a(r)$ some functional inequalities are obtained in \cite{QH}, some interesting monotonicity properties in \cite{WZQC}, and several sharp inequalities in \cite{ZWC2}.

The {\it invariant} of the linear differential equation
\begin{equation}
\label{eq:wpw}
w''+pw'+qw=0,
\end{equation}
where $p$ and $q$ are functions of $z$, is defined to be
$$
I\equiv q- \frac{1}{2}p'- \frac{1}{4}p^2
$$
(cf. [Rai2,p.9]).  If  $w_1$ and $w_2$ are two linearly independent
solutions of (\ref{eq:wpw}), then their quotient
$w\equiv w_2/w_1$ satisfies the differential equation
$$
S_w(z)=2I,
$$
where $S_w$ is the Schwarzian derivative
$$
S_w\equiv \left(\frac{w''}{w'}\right)'-
\frac{1}{2}\left(\frac{w''}{w'}\right)^2
$$
and the primes indicate differentiations (cf. \cite[pp. 18,19]{Rai2}).

From these considerations and the fact that $\mathcal {K}_a(r)$
and $\mathcal {K}'_a(r)$ are linearly independent solutions
of (\ref {eq:kymm2}) (see \cite[(1.11)]{AQVV}), it follows that
$w = \mu_a(r)$ satisfies the differential equation
$$
S_w(r) = \frac{-8a(1-a)}{(r')^2} + \frac{1+6r^2-3r^4}{2r^2(r')^4}.
$$


The \emph{generalized elliptic integral
of the second kind} is given by
\begin{align*}
\mathcal{E}_a &=\mathcal{E}_a(r)\equiv \frac{\pi}{2}
F(a-1,1-a;1;r^2)\\
&=(\sin \pi a)\int_0^{\pi/2}(\tan t)^{1-2a}(1-r^2\sin^2t)^{1-a}\,dt\\
\mathcal{E}_a' &=\mathcal{E}_a'(r)=\mathcal{E}_a(r'),\\
\mathcal{E}_a(0) &=\frac{\pi}{2},\quad \mathcal{E}_a(1)=
\frac{\sin (\pi a)}{2(1-a)}.
\end{align*}
For $a = \frac{1}{2}, ~\mathcal{K}_a$ and $\mathcal{E}_a$
reduce to $\mathcal{K}$
and $\mathcal{E}$, respectively, the usual elliptic
integrals of the first and second kind \cite{BF}, respectively.
Likewise $\mu_{1/2}(r) = \mu(r)$, the modulus of the well-known
Gr\"otzsch ring in the plane \cite{LV}.

\begin{cor}\label{cor:KaEa}
        The generalized elliptic integrals  $\mathcal {K}_a$ and
$\mathcal {E}_a$ satisfy the differential equations
\begin{align}
&r(r')^2\frac{d^2\mathcal{K}_a}{dr^2} + (1-3r^2)
\frac{d\mathcal{K}_a}{dr}-4a(1-a)r\mathcal{K}_a=0,\label{eq:dereka}\\
&r(r')^2\frac{d^2\mathcal{E}_a}{dr^2}+(r')^2\frac{d\mathcal{E}_a}{dr}
+4(1-a)^2r\mathcal{E}_a=0,\label{eq:dertoka}
\end{align}
respectively.
\end{cor}

        {\bf Proof.}  These follow from (\ref{eq:kymm2}).\ \ $\square$
\medskip

For $a = \frac{1}{2}$ these reduce to well-known
differential equations \cite[pp.\ 474-475]{AVV6}, \cite{BF}.

\medskip


\section{Identities of Legendre and Elliott}\label{sect:idenLE}

In geometric function theory the complete elliptic integrals
$\mathcal{K}(r)$ and $\mathcal{E}(r)$ play an
important role. These integrals may be defined, respectively, as
\begin{equation*}
\mathcal{K}(r)=
\textstyle{\frac{\pi}{2}} F(\textstyle{\frac{1}{2}},\textstyle{\frac{1}{2}};
1; r^2),\
\mathcal{E}(r)=
\textstyle{\frac{\pi}{2}} F(\textstyle{\frac{1}{2}},
-\textstyle{\frac{1}{2}};1;r^2),
\end{equation*}
for $-1 < r < 1$. These are $\mathcal{K}_a(r)$ and $\mathcal{E}_a(r)$,
respectively, with $a=\frac{1}{2}$. We also consider the functions
\begin{align*}
\mathcal{K}' &=\mathcal{K}'(r)=\mathcal{K}(r'),\quad 0 < r < 1,\\
\mathcal{K}(0) &=\pi /2,\quad \mathcal{K}(1^-)=+\infty,
\end{align*}
and
\begin{equation*}
\mathcal{E}'=\mathcal{E}'(r)=\mathcal{E}(r'),\quad 0 \leqslant r
\leqslant 1,
\end{equation*}
where $r'=\sqrt{1-r^2}$.
For example, these functions occur in the following quasiconformal counterpart
of the Schwarz Lemma \cite{LV}:

\begin{thm}\label{thm:kin}
For $K\in [1,\infty )$, let  $w$  be a $K$-quasiconformal mapping
of the unit disk $D=\{z:|z|<1\}$ into the unit disk $D'=\{w:|w|<1\}$
with $w(0)=0$.  Then
\begin{equation*}
|w(z)| \leqslant \varphi_K(|z|),
\end{equation*}
where
\begin{equation}\label{eqn:vphi}
\varphi_K(r)\equiv \mu^{-1}\big(\frac{1}{K}\mu(r)\big)\quad \text{and}
\quad \mu(r) \equiv \frac{\pi \mathcal{K}'(r)}{2 \mathcal{K}(r)}.
\end{equation}
This result is sharp in the sense that for each $z\in D$ and
$K \in [1,\infty)$
there is an extremal $K$-quasiconformal mapping that takes the unit
disk $D$ onto the unit disk $D'$ with $w(0)=0$ and
$|w(z)|=\varphi_K(|z|)$ {\rm (}see \cite[p. 63]{LV}{\rm )}.
\end{thm}

It is well known \cite{BF} that the complete elliptic integrals
$\mathcal{K}$ and $\mathcal{E}$ satisfy the Legendre relation
\begin{equation}\label{eq:legendre}
\mathcal{E} \mathcal{K}' + \mathcal{E}' \mathcal{K} - \mathcal{K}
\mathcal{K}' = \frac{\pi}{2}.
\end{equation}
For several proofs of \eqref{eq:legendre}
see \cite{Du}.

In 1904, E. B. Elliott \cite{E} (cf. \cite{AVV5}) obtained the following
generalization of this result.

\begin{thm}\label{thm:Ell} If  $a,b,c \geqslant 0$
and $0 < x < 1$ then
\begin{equation}\label{eq:elliott}
F_1F_2 + F_3F_4 - F_2F_3 =
\frac{\Gamma(a+b+1)\Gamma(b+c+1)}{\Gamma (a+b+c+\frac{3}{2})
\Gamma(b + \frac{1}{2})}.
\end{equation}
where
\begin{align*}
F_1 &= F\bigg({\frac{1}{2}} + a, -\frac{1}{2} - c; 1+a+b; x\bigg),\\
F_2 &= F\bigg(\frac{1}{2} - a, \frac{1}{2} + c; 1+b+c; 1 - x\bigg),\\
F_3 &= F\bigg(\frac{1}{2} + a, \frac{1}{2} - c; 1+a+b;x\bigg),\\
F_4 &= F\bigg(-\frac{1}{2} - a, \frac{1}{2} + c; 1+b + c;1 - x\bigg).
\end{align*}
\end{thm}

Clearly \eqref{eq:legendre} is a special case of \eqref{eq:elliott},
when $a = b = c = 0$ and $x=r^2$.
For a discussion of generalizations of Legendre's
Relation see Karatsuba and Vuorinen
\cite{KV} and Balasubramanian, Ponnusamy, Sunanda Naik, and Vuorinen
\cite{BPSV}.

Elliott proved \eqref{eq:elliott} by a clever change of variables in
multiple integrals.  Another proof was suggested without details in
\cite[p.\ 138]{AAR}, and in \cite{AVV7} we provided the missing details.

\bigskip

The generalized elliptic integrals satisfy the identity
\begin{equation}
\label{eq:legkaea}
\mathcal{E}_a\mathcal{K}_a'+\mathcal{E}_a'\mathcal{K}_a-
\mathcal{K}_a\mathcal{ K}_a'=\frac{\pi\sin(\pi a)}{4(1-a)}.
\end{equation}
This follows from Elliott's formula \eqref{eq:elliott} and
contains the classical relation of Legendre \eqref{eq:legendre}
as a special case. See also Lemma \ref{lem:vv11}.

Finally, we record the following formula of
Kummer \cite[p. 63, Form. 30]{Kum}:
$$
\hspace{-1in}F(a,b;a+b-c+1;1-x)F(a+1,b+1;c+1;x)
$$
$$
\qquad \qquad + \frac{c}{a+b-c+1}F(a,b;c;x)F(a+1,b+1;a+b-c+2;1-x)
$$
$$
\hspace{.25in}  = Dx^{-c}(1-x)^{c-a-b-1},\ \ D= \frac{\Gamma(a+b-c+1) \Gamma(c+1)}{\Gamma(a+1) \Gamma(b+1)}.
$$
This formula, like Elliott's identity, may be rewritten in many
different ways if we use the contiguous relations of Gauss.
Note also the special case $c=a+b-c+1.$

\bigskip

\medskip

\section{Approximation of elliptic integrals and perimeter of ellipse}\label{sect:apprell}

Efficient algorithms for the numerical
evaluation of $\mathcal{K}(r)$ and $\mathcal{E}(r)$
are based on the arithmetic-geometric mean iteration of Gauss.
This fact led to some
close majorant/minorant functions for $\mathcal{K}(r)$ in terms
of mean values in \cite{VV}.  Recently, mean iterations derived from
transformation formulas for the hypergeometric functons have been investigated in \cite{HKM}.

Next, let $a$ and $b$ be the semiaxes of an ellipse with
$a >b$ and eccentricity
$e = \sqrt{a^2 - b^2}/a$, and let $L(a,b)$ denote the perimeter of
the ellipse. Without loss of generality we take $a = 1$.
In 1742, Maclaurin (cf. \cite{AB}) determined that
\begin{equation*}
L(1,b)=4\mathcal{E}(e)=2\pi\cdot {}_2F_1(\textstyle{\frac{1}{2}},
-\textstyle{\frac{1}{2}};1;e^2).
\end{equation*}

In 1883, Muir (cf. \cite{AB}) proposed that $L(1,b)$ could be
approximated by the expression $2\pi [(1+b^{3/2})/2]^{2/3}$.
Since this expression has a close resemblance to the power
mean values studied in \cite{VV}, it is natural to study
the sharpness of this approximation.
Close numerical examination of the error in this
approximation led Vuorinen \cite{V} to conjecture that Muir's
approximation is a lower bound for the perimeter.  Letting
$r=\sqrt{1-b^2}$, Vuorinen asked whether
\begin{equation}\label{eq:vuor}
\frac{2}{\pi}\mathcal{E}(r)
= {}_2F_1 \Big(\tfrac{1}{2},- \tfrac{1}{2}; 1;r^2\Big)
\geqslant \Big( \frac{1+(r')^{3/2}}{2} \Big)^{2/3}
\end{equation}
for all $r\in [0,1]$.

In \cite{BPR1} Barnard and his coauthors proved that
inequality \eqref{eq:vuor} is true.
In fact, they expanded both functions into Maclaurin
series and proved that the differences of the corresponding
coefficients of the two series all have the same sign.

Later, the same authors \cite{BPR2} discovered an upper bound for
$\mathcal{E}$ that complements the lower bound in \eqref{eq:vuor}:
\begin{equation}\label{eq:upper}
\frac{2}{\pi}\mathcal{E}(r)
= {}_2F_1 \Big(\tfrac{1}{2}, - \tfrac{1}{2}; 1;r^2\Big)
\leqslant \Big( \frac{1 +(r')^2}{2} \Big)^{1/2},
\quad 0 \leqslant r \leqslant 1.
\end{equation}
See also \cite{BPS}.

In \cite{BPR2} the authors have considered 13 historical
approximations (by Kepler, Euler, Peano, Muir, Ramanujan,
and others) for the perimeter of an ellipse and determined a
linear ordering among them. Their main tool was the following
Lemma \ref{lem:genhyper} on generalized hypergeometric functions.
These functions are
defined by the formula
\begin{equation*}
_pF_q(a_1,a_2,\cdots ,a_p;b_1,b_2,\cdots ,b_q;z)\equiv
1+\sum_{n=1}^{\infty}\frac{\Pi_{i=1}^p(a_i,n)}{\Pi_{j=1}^q(b_j,n)}
\cdot \frac{z^n}{n!},
\end{equation*}
where  $p$ and $q$ are positive integers and in which no
denominator parameter $b_j$ is permitted to
be zero or a negative integer. When $p=2$ and $q=1$, this reduces to
the usual Gaussian hypergeometric function $F(a,b;c;z)$.

Some of this joint research is discussed in the survey paper \cite{BRT}.

\begin{lem}\label{lem:genhyper}
Suppose $a,b > 0$. Then for any $\epsilon$ satisfying
$\frac{ab}{1 + a + b} < \epsilon < 1$,
\begin{equation*}
_3F_2(-n,a,b; 1 + a + b, 1 + \epsilon - n;1) > 0
\end{equation*}
for all integers  $n \geqslant 1$.
\end{lem}

\begin{OtherAppr}
At the end of the preceding section we pointed out that
upper and lower bounds can be found for $\mathcal{K}(r) $
in terms of mean values.
Another source for the approximation of $\mathcal{K}(r)$ is
based on the asymptotic behavior at the singularity $r=1$,
where $\mathcal{K}(r)$ has logarithmic growth. Some of the
approximations motivated by this aspect will be discussed next.

Anderson, Vamanamurthy, and Vuorinen  \cite{AVV3} approximated
$\mathcal{K}(r)$ by the inverse hyperbolic tangent function
$\arth$, obtaining the inequalities
\begin{equation}\label{eq:arth}
\frac{\pi}{2}\Bigg(\frac{\arth r}{r}\Bigg)^{1/2}
< \mathcal{K}(r)
< \frac{\pi}{2}\, \frac{\arth r}{r},
\end{equation}
for $0<r<1$. Further results were proved by Laforgia and Sismondi
\cite{LS}. K\"uhnau \cite{Ku1} and Qiu \cite{Q1} proved that, for
$0 < r < 1$,
\begin{equation*}
\frac{9}{8+r^2} < \frac{\mathcal{K}(r)}{\log (4/r')}.
\end{equation*}

Qiu and Vamanamurthy \cite{QVa} proved that
\begin{equation*}
\frac{\mathcal{K}(r)}{\log (4/r')}
< 1 + \frac{1}{4}(r')^2\quad \text{for}\ 0 < r < 1.
\end{equation*}
Several inequalities for $\mathcal{K}(r)$ are given in
\cite[Theorem 3.21]{AVV6}.
Later Alzer \cite{A3} showed that
\begin{equation*}
1+\Big(\frac{\pi}{4\log 2}-1\Big)(r')^2
<\frac{\mathcal{K}(r)}{\log (4/r')},
\end{equation*}
for $0<r<1$. He also showed that the constants
$\frac{1}{4}$ and $\pi/( 4\log 2)-1$ in the above
inequalities are best possible. The authoritative NIST
handbook \cite{OLBC} lists some of these inequalities in its Section 19.9.

For further refinements, see \cite[(2.24)]{QVu1} and \cite{Be}.

Alzer and Qiu \cite{AlQ} have written a related manuscript in
which, besides proving many inequalities for complete elliptic
integrals, they have refined \eqref{eq:arth} by proving that
\begin{equation*}
\frac{\pi}{2}\Big(\frac{\arth r}{r}\Big)^{3/4}<\mathcal{K}(r)
< \frac{\pi}{2}\, \frac{\arth r}{r}.
\end{equation*}
They also showed that $3/4$ and $1$ are the best exponents for
$(\arth r)/r$ on the left and right, respectively.  Further estimates for complete elliptic integrals have been obtained in \cite{ABa} and \cite{GQ}.

One of the interesting tools used by the authors of \cite{AlQ}  is the following lemma of
Biernaki and Krzy\.{z} \cite{BK} (for a detailed proof see
\cite{PV1}):
\end{OtherAppr}

\begin{lem}\label{lem:BK}
Let $r_n$ and $s_n$, $n = 1,2,\dots$ be real numbers, and let
the power series $R(x) = \sum_{n=1}^{\infty }r_nx^n$ and $S(x) =
\sum_{n=1}^{\infty}s_nx^n$ be convergent for
$|x|<1$. If $s_n > 0$ for $n=1,2,\ldots$, and if
$r_n/s_n$ is strictly increasing {\rm(}resp. decreasing{\rm)} for
$n=1,2,\ldots$, then the function $R/S$ is strictly increasing
{\rm(}resp. decreasing{\rm)} on $(0,1)$.
\end{lem}

\begin{GenEllRmk}
Recently some new estimates for  $\mathcal{K}(r)$ and $\mathcal{E}(r)$ were obtained in \cite{GQ}.  For the case of generalized elliptic
integrals some inequalities are given in \cite{AQVV}, and further properties are found for them in \cite{HLVV} and \cite{HVV}.
B. C. Carlson has introduced some standard forms for elliptic
integrals involving certain symmetric
integrals. Approximations for these functions can be found in
\cite{CG}.   In \cite{KN} and \cite{OLBC} several inequalities
are obtained for elliptic integrals given in the Carlson form.  In \cite{ABa} S. Andr\'as and \'A. Baricz have compared the generalized elliptic integral  $\mathcal{K}_a(r)$ with certain other zero-balanced hypergeometric functions. In his new book Baricz \cite{Bar3}
investigates various properties of power series and
provides refinements for some of the above results, applying, for example, Lemma \ref{lem:BK}. See also
Zhang, Wang and Chu \cite{ZWC1}. Recently, a variant of Lemma \ref{lem:BK}
for the case when the numerator and denominator are  polynomials of the same
degree, was give in \cite{HVV}.  See also
\cite{KS}.
\end{GenEllRmk}

\section{Hypergeometric series as an analytic function}

For rational triples $(a,b,c)$ the hypergeometric function
often can be expressed in terms of elementary functions.
Long lists with such triples containing hundreds of functions
can be found in \cite{PBM}. For example, the functions
$$ f(z)\equiv z F(1,1;2;z)= -\log(1-z)$$
and
$$g(z)\equiv z F(1,1/2;3/2;z^2)=
 \frac{1}{2} \log \left(\frac{1+z}{1-z}\right)$$
have the property that they both map the unit disk into a strip domain.
Observing that they both correspond to the case $c=a+b$ one may ask
(see \cite{PV1,PV2}) whether there exists $\delta >0$ such that
$zF(a,b;a+b;z)$ and $zF(a,b;a+b;z^2)$ with $a,b \in (0,\delta)$
map into a strip domain.

Membership of hypergeometric functions in some special classes of
univalent functions is studied in \cite{PV1,PV2, PV3, BPV2}.

\section{Generalized modular equations}

The argument $r$ is sometimes called the \emph{modulus} of the elliptic integral $\mathcal{K}(r)$; further, for integer values $p=1,2,\ldots ,$ the equation
\begin{equation}
\frac{\mathcal{K}'(s)}{\mathcal{K}(s)}=p\frac{\mathcal{K}'(r)}{\mathcal{K}(r)},\label{eq:mod}
\end{equation}
with $r,s\in(0,1)$, is called the \emph{modular equation of degree} $p$.  If we use the notation
\begin{equation}
\varphi_K(r) \equiv \mu^{-1}(\mu(K)/K),
\end{equation}
where  $\mu$  is the modulus of the well-known
Gr\"otzsch ring in the plane \cite{LV}, then the solution of (\ref{eq:mod}) for $s$ is given by  $s=\varphi_{1/p}(r).$  We will now discuss some of the numerous modular equations or, more precisely, algebraic consequences of the transcendental equation (\ref{eq:mod}) that were found by Ramanujan.  Our discussion is based on a nice survey of Ramanujan's work in \cite[pp. 4--10]{Be3}.

Ramanujan introduced the convenient notation
$$
\alpha=r^2,\ \ \beta=s^2
$$
for use in connection with (\ref{eq:mod}).  In this notation a third-degree modular equation due to Legendre \cite[p. 105]{BB} takes the form
$$
(\alpha\beta )^{1/4}+((1-\alpha )(1-\beta ))^{1/4} = 1,
$$
with $\alpha = r^2$, $\beta =\varphi_{1/3}(r)^2$.  In the next theorem we list some of Ramanujan's modular equations, in his notation.

\begin{thm}
The function $\varphi_K$ satisfies the following identities for $r\in (0,1)$:
\begin{enumerate}
\item  $(\alpha\beta )^{1/2}+(((1-\alpha )(1-\beta ))^{1/2}+2(16\alpha\beta (1-\alpha )(1-\beta ))^{1/6} = 1\newline
  \text{for}\ \alpha = r^2, \beta = \varphi_{1/5}(r)^2.$

\item $(\alpha\beta )^{1/8}+((1-\alpha )(1-\beta ))^{1/8}=1\\
\text{for}\ \alpha = r^2, \beta=\varphi_{1/7}(r)^2.$

\item $(\alpha (1-\gamma ))^{1/8} +(\gamma (1-\alpha))^{1/8} = 2^{1/3}(\beta (1-\beta ))^{1/24}\newline
\text{for}\ \alpha = r^2, \beta=\varphi_{1/3}(r)^2, \gamma = \varphi_{1/9}(r)^2.$

\item $(\alpha\beta )^{1/8}+((1-\alpha )(1-\beta ))^{1/8}+2^{2/3}(\alpha\beta (1-\alpha )(1-\beta ))^{1/24}=1\newline
\text{for}\ \alpha = r^2, \beta = \varphi_{1/23}(r)^2.$

\item $(\frac{1}{2}(1+\sqrt{\alpha\beta}+\sqrt{(1-\alpha )(1-\beta )}))^{1/2}=(\alpha\beta )^{1/8}+((1-\alpha )(1-\beta ))^{1/8}- \newline
(\alpha\beta (1-\alpha )(1-\beta ))^{1/8}\newline
\text{for}\ \alpha=r^2, \beta=\varphi_{1/7}(r)^2\ \text{or\ for}\ \alpha = \varphi_{1/3}(r)^2, \beta = \varphi_{1/5}(r)^2.$
\end{enumerate}
\end{thm}

All of these identities are from \cite{Be3}: (1) is Entry 13 (i) on p. 280, (2) is Entry 19 (i) on p. 314, (3) is Entry 3 (vi) on p. 352, (4) is Entry 15 (i) on p. 411, and (5) is Entry 21 (i) on p. 435.
\vspace{.05in}

In 1995 B. Berndt, S. Bhargava, and F. Garvan published an important paper \cite{BeBG} in which they studied generalized modular equations and gave proofs for numerous statements concerning these equations made by \linebreak Ramanujan in his unpublished notebooks.  No record of Ramanujan's original proofs has remained.  A \emph{generalized modular equation with signature} $1/a$ \emph{and order} (\emph{or degree}) p is
\begin{equation}\label{eq:modgen}
\frac{F(a,1-a;1;1-s^2)}{F(a,1-a;1;s^2)}=p\frac{F(a,1-a;1;1-r^2)}{F(a,1-a;1;r^2)},\ 0 < r < 1.
\end{equation}
Here $F$ is the Gaussian hypergeometric function defined in (\ref{eq:kolmas}).  The word \emph{generalized} alludes to the fact that the parameter $a\in (0,1)$ is arbitrary.  In the classical case, $a = \frac{1}{2}$ and $p$ is a positive integer.  Modular equations were studied extensively by Ramanujan, see \cite{BeBG}, who also gave numerous algebraic identities for the solutions  $s$ of (\ref{eq:modgen}) for some rational values of $a$ such as $\frac{1}{6}, \frac{1}{4}, \frac{1}{3}.$

To rewrite (\ref{eq:modgen}) in a slightly shorter form, we use the decreasing homeomorphism $\mu_a: (0,1)\to (0,\infty )$ defined by
$$
\mu_a(r)\equiv \frac{\pi}{2\sin (\pi a)}\frac{F(a,1-a;1;1-r^2)}{F(a,1-a;1;r^2)},
$$
for $a\in (0,1).$  We can now rewrite (\ref{eq:modgen}) as
\begin{equation}\label{eq:muasmuar}
\mu_a(s)=p\mu_a(r),\ 0<r<1.
\end{equation}
The solution of (\ref{eq:muasmuar}) is then given by
\begin{equation}\label{eq:varphiKa}
s=\varphi_K^a(r)\equiv\mu_a^{-1}(\mu_a(r)/K),\ \ p = 1/K.
\end{equation}
We call $\varphi_K^a(r)$ \emph{the modular function with signature} $1/a$ \emph{and degree} $p=1/K.$  Monotonicity and convexity properties of $\mathcal{K}_a(r)$, $\mathcal{E}_a(r)$, $\varphi_K^a(r)$, and $\mu_a(r)$ and certain combinations of these special functions are established in \cite{WZC}.

	For the parameter $K=1/p$ with $p$ a small positive integer, the function (\ref{eq:varphiKa}) satisfies several algebraic identities.  The main cases studied in \cite{BeBG} are
$$
	a=\frac{1}{6}, \frac{1}{4}, \frac{1}{3},\ \ p = 2, 3, 5, 7, 11,\ldots .
	$$
	For generalized modular equations we use the Ramanujan notation:
	$$
	\alpha\equiv r^2,\ \ \beta \equiv \varphi^a_{1/p}(r)^2.
$$

We next state a few of the numerous identities \cite{BeBG} satisfied by $\varphi^a_{1/p}$ for various values of the parameters $a$ and $p$.

\begin{thm}\rm{(}\cite[Theorem 7.1]{BeBG}\rm{)}
\emph{If} $\beta$ \emph{has degree} $2$  \emph{in the theory of signature} $3$, \emph{then, with} $a = \frac{1}{3}$, $\alpha = r^2$, $\beta = \varphi^a_{1/2}(r)^2$,
$$
(\alpha\beta )^{1/2}+\{(1-\alpha )(1-\beta )\}^{1/3}=1.
$$
\end{thm}

\begin{thm}\rm{(}\cite[Theorem 7.6]{BeBG}\rm{)}
\emph{If} $\beta$ \emph{has degree} $5$ \emph{then, with} $a=\frac{1}{3}$, $\alpha = r^2$, $\beta =\varphi^a_{1/3}(r)^2$,
$$
(\alpha\beta )^{1/3}+\{(1-\alpha )(1-\beta )\}^{1/3} + 3\{\alpha\beta (1-\alpha )(1-\beta )\}^{1/6} = 1.
$$
\end{thm}

\begin{thm}\rm{(}\cite[Theorem 7.8]{BeBG}\rm{)}
\emph{If} $\beta$ \emph{has degree} $11$ \emph{then, with} $a=\frac{1}{3}$, $\alpha=r^2$, $\beta=\varphi^a_{1/11}(r)^2$,
$$
(\alpha\beta )^{1/3}+\{(1-\alpha )(1-\beta )\}^{1/3} + 6\{\alpha\beta (1-\alpha )(1-\beta )\}^{1/6}
$$
$$\hspace{1in}+3\sqrt{3}\{\alpha\beta (1-\alpha )(1-\beta )\}^{1/12}\{(\alpha\beta )^{1/6}+((1-\alpha )(1-\beta ))^{1/6}\} = 1.
$$

\end{thm}

Such results are surprising, because they provide algebraic identities for the modular function $\varphi^a_K$, which itself is defined in terms of the transcendental function $\mu_a(r)$.  It is an interesting open problem to determine which of the modular equations in \cite{BeBG} can be solved algebraically, explicitly in terms of the modular function.

Because of its geometric significance for  geometric function theory (see (\ref{eqn:vphi}), \cite{LV}, \cite{QVu3}, \cite{Ku2}), it is desirable
to give upper bounds for the function $\varphi_K(r), K>1\,.$ There are many
such bounds in the literature; see the survey \cite{AVV8}.
Much less is known about the function $\varphi^a_K(r), K>1\,,$ but some bounds
can be found in \cite{HVV}, \cite{HLVV}, \cite{WZC},  \cite{WZQC}.

\vskip 1cm

\noindent
ANDERSON: \\
  Department of Mathematics \\
Michigan State University \\
     East Lansing, MI 48824, USA \\
      email: {\tt anderson@math.msu.edu}\\
     FAX: +1-517-432-1562\\ [1mm]

\noindent
VUORINEN:\\
     Department of Mathematics \\
     FIN--00014 \\
     University of Turku, FINLAND\\
     e-mail: ~~{\tt vuorinen@utu.fi}\\
     FAX: +358-2-3336595\\

\end{document}